\documentclass[11pt]{article}
\usepackage{latexsym}
\usepackage[utf8]{inputenc}
\usepackage{amsmath,amssymb,amsthm,enumerate,bm,color}

\newtheorem{thm}{Theorem}[section]

\newtheorem{prop}[thm]{Proposition}
\newtheorem{conj}[thm]{Conjecture}

\theoremstyle{remark}
\newtheorem{rem}[thm]{\bf Remark}

\numberwithin{equation}{section}

\newcommand{\beq}{\begin{equation}}
\newcommand{\eeq}{\qed \end{equation}}

\def\1{\mathbf{1}}
\def\2{\mathbf{2}}

\begin{document}

\begin{Large}
\centerline{\bf Asymptotic Directions for the Zero Sets}
\centerline{\bf of the Components of an Electrical Field}
\centerline{\bf from a Finite Number of Point Charges on the Plane}
\centerline {\bf Part II}
\medskip

\centerline{\bf Tam\'as Erd\'elyi, Joseph Rosenblatt, Rebecca Rosenblatt}
 \end{Large}
\medskip

\begin{abstract}  We study the structure of the zero set of a nontrivial finite point charge electrical field $F = (X,Y)$ in the plane 
$\mathbb R^2$. We establish equations satisfied by the possible directions for the zero sets \{X = 0\} and $\{Y = 0\}$ 
separately, and we show that there are only finitely many possible asymptotic directions for both of these 
zero sets. We suspect that the set of asymptotic directions for \{X = 0\} and the set of asymptotic directions for 
$\{Y = 0\}$ are (essentially) distinct.

\end{abstract}

\section{Introduction}
We study the structure of the zero set of a finite point charge electrical field $F = (X,Y)$ in $\mathbb{R}^2$, where the force vector $F(X,Y)$ is given by 
$$F(x,y) = (X(z,y),Y(x,y)) = \sum_{j=1}^M{\frac{a_j}{r_j^2} \, \frac{(x,y) - (x_j,y_j)}{r_j}}$$
with
$$r_j = \|(x,y) - (x_j,y_j)\|_2= ((x-x_j)^2+(y-y_j)^2)^{1/2}\,.$$
That is,
$$X(x,y) = \sum_{j=1}^M{\frac{a_j(x-x_j)}{((x-x_j)^2 + (y-y_j)^2)^{3/2}}} \qquad \text{and} \qquad \sum_{j=1}^M{\frac{a_j(y-y_j)}{((x-x_j)^2 + (y-y_j)^2)^{3/2}}}\,.$$  
In an earlier paper \cite{ERR} we focused on the ``Special Case" where the point charges are on a line. 
By using a linear transformation this case can be reduced to the case 
$$0 < x_1 < x_2 < \cdots < x_M\,, \qquad y_1 = y_2 = \cdots = y_M = 0$$
Building on results in \cite{GNS}, \cite{K}, \cite{KP}, \cite{Gibson}, \cite{WI}, \cite{WII}, \cite{WB} \cite{WIII}, and \cite{Jan}, 
in \cite{ERR} we gave a fairly complete structural information about the zero sets of $X$ and $Y$ for $F = (X,Y)$ in the Special Case.
A highlight of \cite{ERR} states that in the Special Case the zero set of a nontrivial $F = (X,Y)$ contains at most $9M^24^M$ points,
where $M$ is the number of point charges. on a line.  

Let $0 \neq \alpha \in (-\infty,\infty)$. The line $y = \alpha x$ is called an asymptotic direction for a set $A \subset {\mathbb R}^2$
if there are $(p_m,q_m) \in A$ such that
$$\lim_{|p_m| \rightarrow \infty}{\frac{q_m}{p_m}} = \lim_{|q_m| \rightarrow \infty}{{\frac{q_m}{p_m}}} = \alpha\,.$$
The $x$-axis is called an asymptotic direction for a set $A \subset {\mathbb R}^2$ if there are $(p_m,q_m) \in A$ such that
$$\lim_{|p_m| \rightarrow \infty}{\frac{q_m}{p_m}} = 0\,.$$
The $y$-axis is called an asymptotic direction for a set $A \subset {\mathbb R}^2$ if there are $(p_m,q_m) \in A$ such that
$$\lim_{|q_m| \rightarrow \infty}{\frac{p_m}{q_m}} = 0\,.$$
For the sake of brevity we will use the notation $\beta := 1/\alpha$. Note that the case $\beta = 0$ corresponds to the asymptotic direction given by the $y$-axis.

In this paper we describe the possible asymptotic directions for the zero sets of $X$ and $Y$.  
In the following sections we will derive the zeroing equations that allow us to come close to being able to affirm the conjecture below.
\medskip

\begin{conj} The asymptotic directions of the points in the zero sets $X = 0$ and $Y = 0$, in the Type I and Type II domains, are distinct.
Hence, outside some large disk, there are no common zeros in either domain for $X$ and $Y$.  That is, $F:=(X,Y)$ has no zeros far from the origin.
\end{conj}

\medskip

We call the domains $\{(x,y):|x| < |y|\}$ and $\{(x,y):|y| < |x|\}$ Type I domain and Type II domain, respectively.
In fact, for technical reasons, for a fixed $\delta \in (0,1/2)$, we will consider Type I domains
$\{(x,y):|y| < (1-\delta)|x|\}$ and Type II domains $\{(x,y):|x| < (1-\delta)|y|\}$.  

Suppose $a_1,a_2, \ldots,a_M$ are real numbers which are not all zero. Let $0 \leq L$ be the largest integer for which
$$\sum_{j=1}^M{a_jx_j^ly_j^i} = 0, \qquad i+l < L\,.$$
Such an integer $0 \leq L \leq 2M-1$ exists, otherwise
\begin{equation} \notag
\sum_{j=1}^M{a_jr_j^{2k}} = \sum_{j=1}^M{a_j((x-x_j)^2 + (y-y_j)^2)^k} = \sum_{j=1}^M{a_j(x^2 - 2x_jx + x_j^2 + y^2 - 2y_jy + y_j^2)^k} = 0
\end{equation}
for all $k=0,1,\ldots,M-1$ and $(x,y) \in \mathbb{R}^2$, and by choosing a point $(x,y)$ such that the values $r_j$, $j=1,2,\ldots,M$, are distinct,
the non-vanishing property of the Vandermonde determinants would imply that $a_1=a_2=\cdots=a_M=0$. This number $0 \leq L \leq M-1$ plays a key 
role in our investigation. 

\section{Asymptotic Directions for the Zero Sets of the Components}

\subsection{Type I Domain and $X =0$}

\begin{prop}\label{XTypeI} Let $\delta \in (0,1/2)$ be fixed. If $Y(p_m,q_m) = 0$, where  
$$(p_m,q_m) \in D_{\delta} := \left\{(x,y) \in \mathbb{R}^2: \left| \frac xy \right| \leq 1-\delta\,, \quad \frac{|x|+|x_j|}{|y|-|y_j|} < 1\,, 
\quad |y| > |y_j|, \quad j=1,2,\ldots,M \right\}$$
and
$$\lim_{|q_m| \rightarrow \infty}{\frac{p_m}{q_m}} = \beta\,,$$
then
\begin{equation} \notag
0 = \sum_{l=0}^L \left( \sum_{j=1}^M a_j(-x_j)^ly_j^{L-l} \right) \frac {1}{l!(L-l)!} \frac{d^{L-l}}{d\beta^{L-l}} 
\left( \frac {d^l}{d\beta^l} \left( \frac{\beta}{(1 +\beta^2)^{3/2}} \right)\beta^{L+1} \right)\,.
\end{equation}
There are at most $2L+1=4M-1$ values of $\beta \neq 0$ satisfying the above equation.
\end{prop}

\begin{rem}
We could also interpret the zeroing identify by reversing the order that derivatives are taken. That is, we rescale the equation by $\beta^{L-l+2}$
and interpret the factors 
$$(L-l+2n+2)(L-l-1+2n+2) \cdots (1+2n+2)$$ 
as indicating $L-l$ derivatives have been taken, which would leave terms $\beta^{2n+2}$. But
then rescaling by $\beta^{-1}$ would give the zeroing identity
$$0 = \sum_{l=0}^L \left( \sum_{j=1}^M a_j(-x_j)^ly_j^{L-l} \right) \frac {1}{l!(L-l)!} \frac {d^{l}}{d\beta^{l}} 
\left( \frac{d^{L-l}}{d\beta^{L-l}}\left( \frac{\beta^{L-l+2}}{(1 +\beta^2)^{3/2}} \right) \frac {1}{\beta} \right)\,.$$
Note that this precludes using $\beta = 0$.
\end{rem}

\begin{rem} It is not clear which of these two forms is best to use, and it actually is not so obvious why they give the same zero set, except for $\beta = 0$.
\end{rem}

\begin{proof}
We can consider a Type I domain and consider the Taylor expansion for the zeroing equation for $X$. 
As with previous series expansions, we need to use both the binomial and geometric series. By the Binomial Theorem we have 
$$\frac 1{(1+t)^{3/2}} = (1-t)^{-3/2} = \sum_{n=0}^\infty \binom{-3/2}{n} t^n\,, \qquad t \in (-1,1)\,.$$
Hence
\begin{equation} \notag \begin{split}
\pm X(x,y) & = \, \sum_{j=1}^M{\frac{a_j(x-x_j)}{(y-y_j)^3} \left( 1 + \left( \displaystyle{ \frac{x-x_j}{y-y_j}} \right)^2 \right)^{-3/2}} = 
\sum_{n=0}^\infty{\binom{-3/2}{n} \sum_{j=1}^M{\frac{a_j(x-x_j)^{2n+1}}{(y-y_j)^{2n+3}}}} \cr 
& = \, \sum_{n=0}^\infty{\binom{-3/2}{n} \frac{1}{y^{2n+3}} \sum_{j=1}^M{\frac{a_j(x-x_j)^{2n+1}}{(1-y_j/y)^{2n+3}}}}\,, \qquad \left| \frac{x-x_j}{y-y_j} \right| < 1\,, \quad |y| > |y_j|\,.\cr
\end{split} \end{equation}
This equation is like the one for $X$ from the Special Case analysis in \cite{ERR}. Indeed, now both equations $Y = 0$ and $X = 0$ have the same combined nature,
a need for a binomial expansion for the numerator terms $(x-x_j)^{2n}$ and a geometric powers series expansion for $(1-y_j/y)^{-(2n+3)}$.
So we carry out these expansions.

The Binomial Theorem gives that
$$\frac{1}{(1-y_j/y)^{2n+3}} = \sum_{i=0}^\infty{\binom{-(2n+3)}{i}\left( \frac{-y_j}{y} \right)^i} = 
\sum_{i=0}^\infty \binom{i+2n+2}{2n+2}\left (\frac {y_j}y\right )^i\,,  \qquad |y| > |y_j|\,.$$
So we have
\begin{equation} \label{F1} \begin{split}
X(x,y) & = \, \sum_{n=0}^\infty{\binom{-3/2}{n} \frac{1}{y^{2n+3}} \sum_{j=1}^M{a_j\sum_{l=0}^{2n+1}\binom{2n+1}{l}(-x_j)^lx^{2n+1-l}\sum_{i=0}^\infty \binom{i+2n+2}{2n+2}\left (\frac {y_j}y\right )^i}} \cr
& = \, \sum_{n=0}^\infty{\binom{-3/2}{n} \sum_{l=0}^{2n+1}\sum_{i=0}^\infty \binom{2n+1}l \binom{i+2n+2}{2n+2} \left (\sum_{j=1}^M a_j(-x_j)^ly_j^i \right ) \frac {x^{2n+1-l}}{y^{2n+3 +i}}} \cr 
\end{split} \end{equation}
for all $(x,y) \in D_\delta$. Now for our asymptotic analysis, we need to identify values of $l$ and $i$ such that the moment $\sum_{j=1}^M a_jx_j^ly_j^i$ is zero and not zero.
To eliminate the terms that are zero, we take the moments with smallest $L = l+i$ for which at least one of these moments is not zero, 
say with $l= l_0$ and $i=i_0$. We have  
$$\sum_{n=0}^\infty{\binom{-3/2}{n} \sum_{l=0}^{2n+1}\sum_{i=0}^\infty \binom{2n+1}l \binom{i+2n+2}{2n+2} \left(\sum_{j=1}^M a_j(-x_j)^ly_j^i\right) 
\frac{x^{2n+1}}{y^{2n+1}}\frac {x^{-l}}{y^{-l}}\frac 1{y^{l+i+2}}} = 0$$
for all $(x,y) \in D_\delta$, where the sums over $l$ and $i$ are over all pairs $(l,i)$ with $l+i \geq L = l_0+i_0$.
The above triple sum converges absolutely in $D_\delta$ as
\begin{equation} \notag \begin{split}
& \, \sum_{n=0}^\infty{\left| \binom{-3/2}{n} \right| \sum_{l=0}^{2n+1}\sum_{i=0}^\infty \binom{2n+1}l \binom{i+2n+2}{2n+2} \left|\sum_{j=1}^M a_j(-x_j)^ly_j^i\right| 
\left| \frac{x^{2n+1}}{y^{2n+1}}\frac {x^{-l}}{y^{-l}}\frac 1{y^{l+i+2}}\right|} \cr
\leq & \, \sum_{n=0}^\infty \left| \binom{-3/2}{n} \right| \sum_{l=0}^{2n+1} \sum_{i=0}^\infty \binom{2n+1}l \binom{i+2n+2}{2n+2} \left( \sum_{j=1}^M{|a_j||x_j|^l|y_j|^i} \right)  
\left| \frac{x^{2n+1}}{y^{2n+1}} \frac{x^{-l}}{y^{-l}} \frac{1}{y^{l+i+2}} \right| \cr
= & \, \sum_{j=1}^M{\frac{|a_j|(|x|+|x_j|)}{\left( (|x|+|x_j)^2 + (|y| - |y_j|)^2 \right)^{3/2}}}\,, \qquad (x,y) \in D_\delta\,. \cr
\end{split} \end{equation}
Now assume that $(p_m,q_m) \in D_{\delta}$, $X(p_m,q_m) = 0$, $|q_m| \geq 1$, $m=1,2,\ldots$,  
$$\lim_{m \rightarrow \infty}{\frac{p_m}{q_m}} = \beta \in [-1+\delta,1-\delta] \qquad \text{and} \qquad \lim_{m \rightarrow \infty}{|q_m|} = \infty\,.$$ 
We plug $(p_m,q_m)$ in \ref{F1}, multiply by $q_m^{L+2}$, and separate the terms in which $l+i=k=L$, in which $L+1 \leq l+i = k \leq 2n+1$, and in which $l+i = k \geq N := \max\{2n+2,L+1\}$ 
to obtain 
\begin{equation} \label{F2} \begin{split}
0 & = \, q_m^{L+2} X(p_m,q_m) \cr 
& = \, \sum_{n=0}^\infty \binom{-3/2}{n} \sum_{l=0}^{2n+1} \sum_{i=0}^\infty \binom{2n+1}{l} \binom{i+2n+2}{2n+2} \left( \sum_{j=1}^M{a_j(-x_j)^ly_j^i} \right) 
\left( \frac{p_m}{q_m} \right)^{2n+1-l} \frac{q_m^{L+2}}{q_m^{l+i+2}} \cr
& = \, \sum_{n=0}^\infty{\binom{-3/2}{n} \sum_{l=0}^L \binom{2n+1}{l} \binom{L-l+2n+2}{2n+2} \left( \sum_{j=1}^M{a_j(-x_j)^ly_j^{L-l}} \right) \left( \frac{p_m}{q_m} \right)^{2n+1-l}} \cr 
& + \, \sum_{n=0}^\infty{\binom{-3/2}{n} \sum_{k=L+1}^{2n+1}\sum_{l=0}^k \binom{2n+1}{l} \binom{k-l+2n+2}{2n+2} 
\left(\sum_{j=1}^M{a_j(-x_j)^{l}y_j^{k-l}} \right) \left( \frac{p_m}{q_m} \right)^{2n+1-l}} \frac{1}{q_m^{k-L}} \cr
& + \, \sum_{n=0}^\infty{\binom{-3/2}{n} \sum_{k=N}^{\infty}\sum_{l=0}^{2n+1} \binom{2n+1}{l} \binom{k-l+2n+2}{2n+2} 
\left(\sum_{j=1}^M{a_j(-x_j)^{l}y_j^{k-l}} \right) \left( \frac{p_m}{q_m} \right)^{2n+1-l}} \frac{1}{q_m^{k-L}} \cr
\end{split} \end{equation}
We now show that the last two sums in \ref{F2} multiplied by $|q_m|^{1/2}$ is uniformly bounded when $(p_m,q_m) \in D_{\delta}$ and $|q_m| \geq c$ 
with a constant $c$ depending only on $\delta$, 
$$\gamma := 1 + \max_{1 \leq j \leq M}{(|x_j| + |y_j|)} \qquad \text{and} \qquad A := \max_{1 \leq j \leq M}{|a_j|}\,.$$
We will need the estimate
\begin{equation} \label{F3} \begin{split}
\sum_{n=0}^\infty{(n+1)^kt^n} \leq \sum_{n=0}^\infty{(n+1)(n+2) \cdots (n+k)t^n} = \frac{k!}{(1-t)^{k+1}}, \qquad t \in [0,1)\,.
\end{split} \end{equation}
Observe that $\displaystyle{\frac{k-L-1/2}{2L+4} \geq k+1}$ for $k \geq L+1$, and hence with the notation $s_m := |q_m|^{\frac{1}{2L+4}}$ we have 
\begin{equation} \label{F4} \begin{split}
& \left| \sum_{n=0}^\infty \binom{-3/2}{n} \sum_{k=L+1}^{2n+1} \sum_{l=0}^k \binom{2n+1}{l} \binom{k-l+2n+2}{2n+2} 
\left(\sum_{j=1}^M{a_j(-x_j)^{l}y_j^{k-l}} \right) \left( \frac{p_m}{q_m} \right)^{2n+1-l} \frac{|q_m|^{1/2}}{q_m^{k-L}} \right| \cr
& \quad \leq \sum_{n=0}^\infty (n+1) \sum_{k=L+1}^{2n+1} \sum_{l=0}^k{\frac{2n+1)^{l}}{l!} \frac{(k-l+2n+2)^{k-l}}{(k-l)!} A\gamma^k (1-\delta)^{2n+1-l}} \frac{1}{|q_m|^{k-L-1/2}} \cr
& \quad \leq \sum_{n=0}^\infty (n+1) \sum_{k=L+1}^{2n+1} \sum_{l=0}^k{\frac{(4n+3)^k}{l!(k-l)!} A\gamma^k (1-\delta)^{2n+1-l}} \frac{1}{s_m^{k+1}} \cr
& \quad \leq \sum_{n=0}^\infty (n+1) \sum_{k=L+1}^{2n+1} \sum_{l=0}^k{\frac{1}{k!} \binom{k}{l} (4n+3)^k A\gamma^k (1-\delta)^{2n+1-l}} \frac{1}{s_m^{k+1}} \cr
& \quad \leq \sum_{n=0}^\infty (n+1) \sum_{k=L+1}^{2n+1} \frac{2^k}{k!} (4n+3)^k A\gamma^k (1-\delta)^{2n+1-l} \frac{1}{s_m^{k+1}}\,. \cr
\end{split} \end{equation}
Now we break the sum in the last line of \ref{F3} for $L+1 \leq k \leq n$ and $n+1 \leq k \leq 2n+1$. Recalling \ref{F3} we have  
\begin{equation} \label{F5} \begin{split}
& \, \sum_{n=0}^\infty (n+1) \sum_{k=L+1}^{n} \frac{2^k}{k!} (4n+3)^k A\gamma^k (1-\delta)^{2n+1-l} \frac{1}{s_m^{k+1}} \cr 
\leq & \sum_{n=0}^\infty \sum_{k=L+1}^{n} \frac{2^k}{k!} (4(n+1))^{k+2} A\gamma^k (1-\delta)^{n+1} \frac{1}{s_m^{k+1}} \cr
\leq & \sum_{n=0}^\infty \sum_{k=L+1}^{\infty} A \frac{(8\gamma(n+1))^{k+2}}{k!} (1-\delta)^{n+1} \frac{1}{s_m^{k+1}} \cr  
= & \sum_{k=L+1}^\infty \sum_{n=0}^{\infty} A \frac{(8\gamma(n+1))^{k+1}}{k!} (1-\delta)^{n+1} \frac{1}{s_m^{k+1}} \cr
= & \sum_{k=L+1}^\infty \sum_{n=0}^{\infty} \frac{A}{k!} \left( \frac{8\gamma(n+1)}{s_m} \right)^{k+2} (1-\delta)^{n+1} \cr
\leq & \sum_{k=L+1}^\infty \sum_{n=0}^{\infty} \frac{A}{(k+1)!} (n+1)^{k+1} (1-\delta)^{n+1} \left( \frac{16\gamma}{s_m} \right)^{k+1} \cr
\leq & \sum_{k=L+1}^\infty A \left( \frac{1}{\delta} \right)^{k+2} \left( \frac{16\gamma}{s_m} \right)^{k+2} \cr
\leq & \sum_{k=L+1}^\infty A \left( \frac{16\gamma}{s_m\delta} \right)^{k+2} \leq B_1\,,\cr
\end{split} \end{equation}
with a constant $B_1$ depending only on $\delta$, $A$, and $\gamma$ if $s_m := |q_m|^{\frac{1}{2L+4}} \geq 32\gamma/\delta$. 
 
Similarly
\begin{equation} \label{F6} \begin{split}
& \, \sum_{n=0}^\infty (n+1) \sum_{k=n+1}^{2n+1} \frac{2^k}{k!} (4n+3)^k A\gamma^k (1-\delta)^{2n+1-l} \frac{1}{s_m^{k+1}} \cr 
\leq & \sum_{n=0}^\infty \sum_{k=n+1}^{2n+1} \frac{2^k}{k!} (4(n+1))^{k+1} A\gamma^k \frac{1}{s_m^{k+1}} 
\leq \sum_{n=0}^\infty \sum_{k=n+1}^{2n+1} A \frac{(8\gamma(n+1))^{k+2}}{k!} \frac{1}{s_m^{k+2}} \cr  
\leq & \sum_{n=0}^\infty \sum_{k=n+1}^{2n+1} A \left( \frac{16e\gamma(n+1)}{ks_m} \right)^{k+1} 
\leq \sum_{n=0}^\infty \sum_{k=n+1}^{2n+1} A \left( \frac{16e\gamma(n+1)}{{(n+1)s_m}} \right)^{k+1} \cr 
\leq & \sum_{n=0}^\infty \sum_{k=n+1}^{2n+1} A \left( \frac{32e\gamma}{s_m} \right)^{k+1}2^{-n} \leq B_2 \cr
\end{split} \end{equation}
with a constant $B_2$ depending only on $\delta$, $A$, and $\gamma$ if $s_m := |q_m|^{\frac{1}{2L+4}} \geq 64e\gamma$.
Combining \ref{F4}, \ref{F5}, and \ref{F6} completes the proof that the last but one sum  
\begin{equation} \notag 
\sum_{n=0}^\infty{\binom{-3/2}{n} \sum_{k=L+1}^{2n+1}\sum_{l=0}^k \binom{2n+1}{l} \binom{k-l+2n+2}{2n+2} 
\left(\sum_{j=1}^M{a_j(-x_j)^{l}y_j^{k-l}} \right) \left( \frac{p_m}{q_m} \right)^{2n+1-l}} \frac{1}{q_m^{k-L}} 
\end{equation}
in \ref{F2} converges to $0$ when $\displaystyle{\lim_{m rightarrow \infty}{|q_m|}} = \infty$.

With $N := \max\{2n+2,L+1\}$ and $s_m := |q_m|^{\frac{1}{2L+4}}$ we also have
\begin{equation} \label{F7} \begin{split}
& \, \left| \sum_{n=0}^\infty{\binom{-3/2}{n} \sum_{k=N}^\infty \sum_{l=0}^{2n+1} \binom{2n+1}{l} \binom{k-l+2n+2}{2n+2} 
\left(\sum_{j=1}^M{a_j(-x_j)^{l}y_j^{k-l}} \right) \left( \frac{p_m}{q_m} \right)^{2n+1-l}} \frac{|q_m|^{1/2}}{q_m^{k-L}} \right| \cr
\leq & \,\sum_{n=0}^\infty (n+1) \sum_{k=N}^\infty \sum_{l=0}^{2n+1} \binom{2n+1}{l} \binom{k-l+2n+2}{2n+2}  A\gamma^k \frac{1}{|q_m|^{k-L-1/2}} \cr
\leq & \, \sum_{n=0}^\infty (n+1) \sum_{k=N}^{\infty} \sum_{l=0}^{2n+1} \binom{2n+1}{l} \frac{(k+2n+2)^{2n+2}}{(2n+2)!} A\gamma^k \frac{1}{s_m^{k+1}} \cr
\leq & \, \sum_{n=0}^\infty (n+1) \sum_{k=N}^{\infty} 2^{2n+1} \frac{(k+2n+2)^{2n+2}}{(2n+2)!} A\gamma^k  \frac{1}{s_m^{k+1}} \cr
\leq & \, \sum_{n=0}^\infty \sum_{k=N}^{\infty} (n+1) \left( \frac{2e(2n+2+k)}{2n+2} \right)^{2n+2} A\gamma^k  \frac{1}{s_m^{k+1}} \cr
\leq & \, \sum_{n=0}^\infty \sum_{k=N}^{\infty} (n+1) (2e)^{2n+2}\left( 1 + \frac{k}{2n+2} \right)^{2n+2} A\gamma^k  \frac{1}{s_m^{k+1}} \cr
\leq & \, \sum_{n=0}^\infty \sum_{k=N}^{\infty} (n+1) (2e)^{2n+2} e^k A\gamma^k  \frac{1}{s_m^{k+1}} 
\leq \sum_{n=0}^\infty \sum_{k=N}^{\infty} (n+1) (2e)^{2n+2} A \left( \frac{e\gamma}{s_m} \right)^k \frac{1}{s_m} \cr  
\leq & \, \sum_{n=0}^\infty \sum_{k=N}^{\infty} (n+1) A \left( \frac{4e^2\gamma}{s_m} \right)^k (2e)^{-(2n+2)}\frac{1}{s_m} \leq B_3 \cr
\end{split} \end{equation}
with a constant $B_3$ depending only on $\delta$, $A$, and $\gamma$ if $s_m := |q_m|^{\frac{1}{2L+4}} \geq 8e^2\gamma$. It follows from \ref{F7} that the last sum 
\begin{equation} \notag
\sum_{n=0}^\infty{\binom{-3/2}{n} \sum_{k=2n+2}^{\infty} \sum_{l=0}^{2n+1} \binom{2n+1}{l} \binom{k-l+2n+2}{2n+2} 
\left(\sum_{j=1}^M{a_j(-x_j)^{l}y_j^{k-l}} \right) \left( \frac{p_m}{q_m} \right)^{2n+1-l}} \frac{1}{q_m^{k-L}} 
\end{equation}
in \ref{F2} converges to $0$ when $q_m$ tends to $\infty$.

Now we are ready to take the limit in \ref{F2} when $\displaystyle{\lim_{m \rightarrow \infty}{|q_m|} = \infty}$ and $\displaystyle{\lim_{m \rightarrow \infty}{p_m/q_m} = \beta\,.}$ 
We have seen that the last two sums in \ref{F2} tend to $0$ as $\displaystyle{\lim_{m \rightarrow \infty}{|q_m|} = \infty}$. So we have 
\begin{equation} \notag \begin{split}
0 & = \, \lim_{|q_m| \rightarrow \infty} \sum_{n=0}^\infty{\binom{-3/2}{n} \sum_{l=0}^L \binom{2n+1}{l} \binom{L-l+2n+2}{2n+2} \left( \sum_{j=1}^M{a_j(-x_j)^ly_j^{L-l}} \right) 
\left( \frac{p_m}{q_m} \right)^{2n+1-l}} \cr  
& = \, \sum_{l=0}^L \left( \sum_{j=1}^M{a_j(-x_j)^ly_j^{L-l}} \right) \beta^{-l} \sum_{n=0}^\infty \binom{-3/2}{n} \binom{2n+1}{l} \binom{L-l+2n+2}{2n+2} \beta^{2n+1}\,. \cr
\end{split} \end{equation}
Observe that the right-hand side is a not identically zero Laurent series in $0 \neq \beta \in (-1,1)$.  
We would like to show this zeroing formula is equivalent to a formula involving derivatives in such a manner that it can be seen that there are only a finite number of solutions.
We can do this by observing that
\begin{equation} \notag  
\sum_{n=0}^\infty{\binom{-3/2}{n} \binom{2n+1}{l} \binom{L-l+2n+2}{2n+2} \beta^{2n+1}}
\end{equation}
\begin{equation} \notag 
= \frac {1}{l!(L-l)!} \sum_{n=0}^\infty \binom{-3/2}{n} (2n+1) \cdots (2n +1 -l+1) \cdot (L-l+2n+2)\cdots (1+2n+2) \beta^{2n+1}\,.
\end{equation} 
Interpret the factors 
$$(2n+1)(2n+1-1) \cdots (2n+1-l+1)$$ 
as indicating that the $l$th derivatives of $\displaystyle{\frac \beta{(1+\beta^2)^{3/2}}}$ have been taken, which would leave terms with $\beta^{2n - l +1}$.
But then rescaling by $\beta^{L+1}$, would give this series to be
$$\frac 1{\beta} \frac {d^{L-l}}{d\beta^{L-l}} \left (\frac {d^l}{d\beta^l}\left (\frac {\beta}{(1 +\beta^2)^{3/2}}\right )\beta^{L+1}\right )\,.$$
So our zeroing identity is
\begin{equation} \notag
\sum_{l=0}^L \left( \sum_{j=1}^M a_j(-x_j)^ly_j^{L-l} \right) \frac {1}{l!(L-l)!} \frac{d^{L-l}}{d\beta^{L-l}} 
\left( \frac {d^l}{d\beta^l} \left( \frac{\beta}{(1 +\beta^2)^{3/2}} \right)\beta^{L+1} \right) = 0\,.
\end{equation}
This shows that the zeroing equation is of the form $P(\beta)/(1 +\beta^2)^{(2L+3)/2}$ for a suitable not identically zero polynomial $P$ of degree at most $2L+1 \leq 4M-1$.
So the zeroing identity only gives zeros in the roots of $P$. As $P$ is a polynomial of degree at most $2L+1 \leq 4M-1$, there are at most $2L+1 \leq 4M-1$ 
possible asymptotic directions for $X = 0$ in the Type I domain.
\end{proof}

\medskip

\subsection{Type I Domain and $Y = 0$}

There is more symmetry in our arguments with charges for the General Case. When we switch from $X = 0$ to $Y =0$ in the Type I domain, it just decreases the power in 
the binomial terms by $1$ and decreases the powers in the geometric terms by $1$. With the necessary adjustment of the powers in the Type I domain for $Y = 0$, 
in the Type I domain for $X = 0$ we have the zeroing identity 
$$0 = \sum_{l=0}^L \left( \sum_{j=1}^M a_j(-x_j)^ly_j^{L-l} \right) \beta^{-l} \sum_{n=0}^\infty \binom{-3/2}{n} \binom{2n}l \binom{L-l+2n+1}{2n+1} \beta^{2n}\,,$$
and the following result holds. 

\begin{prop}\label{YTypeI} Let $\delta \in (0,1/2)$ be fixed. If $X(p_m,q_m) = 0$, where
$$(p_m,q_m) \in D_{\delta} := \left\{(x,y) \in \mathbb{R}^2: \left| \frac xy \right| \leq 1-\delta\,, \quad \frac{|x|+|x_j|}{|y|-|y_j|} < 1\,, 
\quad |y| > |y_j|, \quad j=1,2,\ldots,M \right\}$$
and
$$\lim_{|q_m| \rightarrow \infty}{\frac{p_m}{q_m}} = \beta\,,$$
then
\begin{equation} \notag
0 = \sum_{l=0}^L \left( \sum_{j=1}^M a_j(-x_j)^ly_j^{L-l} \right) \frac{1}{l!(L-l)!} 
\frac{d^{L-l}}{d\beta^{L-l}} \left( \frac{d^l}{d\beta^l} \left( \frac{1}{(1 +\beta^2)^{3/2}} \right) \beta^{L+1} \right)\,.
\end{equation}
There are at most $2L+1=4M-1$ values of $\beta \neq 0$ satisfying the above equation.
\end{prop}

\subsection{Common Asymptotic Directions of the Components $X$ and $Y$ in Type I Domain}

We speculate that there are no $\beta$ which is a common asymptotic direction for both $\{X=0\}$ and $\{Y=0\}$ in the Type I domain.    
Such a $\beta \in [-1+\delta,\delta]$ must satisfy both of the equations 
$$0 = \sum_{l=0}^L \Gamma(l) \frac{d^{L-l}}{d\beta^{L-l}} \left( \frac{d^l}{d\beta^l} \left( \frac{\beta}{(1 +\beta^2)^{3/2}} \right) \beta^{L+1} \right)$$ 
and
$$0 = \sum_{l=0}^L \Gamma(l) \frac{d^{L-l}}{d\beta^{L-l}} \left( \frac{d^l}{d\beta^l} \left( \frac{1}{(1 +\beta^2)^{3/2}} \right) \beta^{L+1} \right)\,.$$ 
with some real numbers $\Gamma(l)$ which are not all zero.
In \cite{ERR} we were able to show this in the ``Special Case" where the point charges are on the positive $x$-axis. In fact, in \cite{ERR} 
we showed that in the ``Special Case" the possible asymptotic directions for $\{X=0\}$ and the possible asymptotic directions for$\{Y=0\}$ are 
strictly interlacing in the Type I domain.     

Examples show that the possible asymptotic directions of $X = 0$ and $Y = 0$ in the Type I domain are not necessarily strictly interlacing, 
but we conjecture that they are distinct directions, and hence, outside some large disk $X$ and $Y$ have no common zeros in the Type I domain,  
that is, $F = (X,Y)$ has no zeros far from the origin in the Type I domain.

\medskip

\subsection{Type II Domain and $X = 0$}

The result in this section can be obtained from the result in \ref{XTypeI} by replacing the roles of $x$ and $y$ and $x_j$ and $y_j$.

\begin{prop}\label{XTypeII} Let $\delta \in (0,1/2)$ be fixed. If $X(p_m,q_m) = 0$, where
$$(p_m,q_m) \in H_{\delta} := \left\{(x,y) \in \mathbb{R}^2: \left| \frac yx \right| \leq 1-\delta\,, \quad \frac{|y|+|y_j|}{|x|-|x_j|} < 1\,, 
\quad |x| > |x_j|\,, \quad j=1,2,\ldots,M \right\}$$
and
$$\lim_{|p_m| \rightarrow \infty}{\frac{q_m}{p_m}} = \alpha\,,$$
then
\begin{equation} \notag
0 = \sum_{l=0}^L \left( \sum_{j=1}^M a_j(-x_j)^{L-l}y_j^{l} \right) \frac {1}{l!(L-l)!} \frac{d^l}{d\beta^l} 
\left( \frac{d^{L-l}}{d\alpha^{L-l}} \left( \frac{\alpha}{(1 +\alpha^2)^{3/2}} \right)\alpha^{L+1} \right)\,.
\end{equation}
There are at most $2L+1=4M-1$ values of $\alpha \neq 0$ satisfying the above equation.
\end{prop}

\subsection{Type II Domain and $Y = 0$}

The result in this section can be obtained from the result in \ref{YTypeI} by replacing the roles of $x$ and $y$ and $x_j$ and $y_j$.

\begin{prop}\label{YTypeII} Let $\delta \in (0,1/2)$ be fixed. If $X(p_m,q_m) = 0$, where
$$(p_m,q_m) \in H_{\delta} := \left\{(x,y) \in \mathbb{R}^2: \left| \frac yx \right| \leq 1-\delta\,, \quad \frac{|y|+|y_j|}{|x|-|x_j|} < 1\,, 
\quad |x| > |x_j|\,, \quad j=1,2,\ldots,M \right\}$$
and
$$\lim_{|p_m| \rightarrow \infty}{\frac{q_m}{p_m}} = \alpha\,,$$
then
\begin{equation} \notag
0 = \sum_{l=0}^L \left( \sum_{j=1}^M a_j(-x_j)^{L-l}y_j^{l} \right) \frac {1}{l!(L-l)!} \frac{d^l}{d\beta^l} 
\left( \frac{d^{L-l}}{d\alpha^{L-l}} \left( \frac{1}{(1 +\alpha^2)^{3/2}} \right)\alpha^{L+1} \right)\,.
\end{equation}
There are at most $2L+1=4M-1$ values of $\alpha \neq 0$ satisfying the above equation.
\end{prop}

\subsection{Common Asymptotic Directions of the Components $X$ and $Y$ in Type II domain}

We speculate that there are no $\alpha$ which is a common asymptotic direction for both $\{X=0\}$ and $\{Y=0\}$ in the Type II domain.
Such an $\alpha \in [-1+\delta,\delta]$ must satisfy both of the equations
$$0 = \sum_{l=0}^L \Delta(l) \frac {d^l}{d\alpha^l} \left( \frac{d^{L-l}}{d\alpha^{L-l}} \left( \frac{\alpha}{(1 +\alpha^2)^{3/2}} \right) \alpha^{L+1} \right)$$
and
$$0 = \sum_{l=0}^L \Delta(l) \frac{d^l}{d\alpha^l} \left( \frac{d^{L-l}}{d\alpha^{L-l}}\left( \frac{1}{(1 +\alpha^2)^{3/2}} \right) \alpha^{L+1} \right)$$
with some real numbers $\Delta(l)$ which are not all zero.

In \cite{ERR} we were able to show this in the ``Special Case" where the point charges are on the positive $x$-axis. In fact, in \cite{ERR}    
we showed that in the ``Special Case" the possible asymptotic directions for $\{X=0\}$ and the possible asymptotic directions for$\{Y=0\}$ are 
strictly interlacing in the Type II domain.

Examples show that the possible asymptotic directions of $X = 0$ and $Y = 0$ in the Type II domain are not necessarily strictly interlacing,  
but we conjecture that they are distinct directions, and hence, outside some large disk $X$ and $Y$ have no common zeros in the Type II domain, 
that is, $F = (X,Y)$ has no zeros far from the origin in the Type II domain.

\medskip

\subsection{A Question about Common Zeros of Certain Polynomials}\label{zerosummary}

Let us summarize the basic challenge of distinguishing the asymptotic directions of the zero sets of $X$ and $Y$ in the General Case. 
Let $D^jF = d^jF/dx^j$ denote the usual $j$-th derivative of $F$ with respect to $F$. Take $s \ge 2$. Consider a nontrivial linear combination 
$$\Sigma(R) = \sum_{l=0}^L{c_l D^l[x^{L+1} D^{L-l}(R)]}\,.$$ 
where $c_l$ are real numbers which are not all zero.

\begin{conj}\label{Conjecture1} Let  
$$R_1(x) := A_0(x) := \frac{x}{(1+x^2)^{3/2}} \qquad \text{and} \qquad R_2(x) := B_0(x) := \frac{1}{(1+x^2)^{3/2}}\,.$$
The positive zeros of $\Sigma(R_1)$ and $\Sigma(R_2)$ are distinct.
\end{conj}

\medskip

\noindent Observe that $\Sigma(R_1)$ is an odd function and $\Sigma(R_2)$ is an even function. So the above conjecture would imply that 
the negative zeros of $\Sigma(R_1)$ and $\Sigma(R_2)$ are distinct. Both $\Sigma(R_1)$ and $\Sigma(R_2)$ have zeros of some order at $0$.  

\medskip

\noindent Here are some computations and wishful thinking that suggest that this conjecture is true.
Suppose that either $R := R_1$ or $R := R_2$. We can rewrite $\Sigma(R)$ as follows. First, for $0 \leq i \leq s$,
$$D^l[x^{L+1} D^{L-l}(R)]  = \sum_{k=0}^l \binom lk D^k(x^{L+1})D^{l-k}(D^{L-l}(R)) = \sum_{k=0}^l \binom lk D^k(x^{L+1})D^{L-k}(R)\,.$$
So we have
\begin{equation} \notag \begin{split}
\Sigma(R) & = \, \sum_{l=0}^s{c_l \sum_{k=0}^l{\binom lk D^k(x^{L+1})D^{L-k}(R)}} = \sum_{k=0}^s{\sum_{l=k}^s{\binom lk c_lD^k(x^{L+1})D^{L-k}(R)}} \cr 
& = \, \sum_{k=0}^L{\sigma_kD^k(x^{L+1})D^{L-k}(R)}\,, \cr 
\end{split} \end{equation}
where 
$$\sigma_k := \sum_{l=k}^s{\binom lk c_l}\,.$$  
So 
$$\Sigma(R)(x) = \sum_{k=0}^L{d_k x^{L+1-k}D^{L-k}(R)(x)}\,,$$
where $d_0 := \sigma_0$ and $d_k := \sigma_k \cdot (L+1) \cdots (L+1- k+1)$ for $1 \leq k \leq L$.  
Hence, 
$$\Sigma(R)(x) = \sum_{l=0}^L{\gamma_l x^{l+1}D^l(R)(x)}\,,$$
where for simplicity we have replaced $d_{L-l}$ by $\gamma_l$.
It seems reasonable to use this simpler form for $\Sigma(R)$ to prove the conjecture. Because we are taking $c_l$ to be arbitrary real numbers, 
we must assume here that $\gamma_l$ are arbitrary real numbers too. Again, consider 
$$R_1(x) := A_0(x) := \frac{x}{(1+x^2)^{3/2}} \qquad \text{and} \qquad R_2(x) := B_0(x) := \frac{1}{(1+x^2)^{3/2}}\,.$$  
Let $A_l$ and $B_l$ be the $l$-th derivatives of $A_0$ and $B_0$, respectively, that is,  
$$A_l(x) := \frac{d^lA_0}{dx^l} \qquad \text{and} \qquad B_l(x) := \frac{d^lB_0}{dx^l}\,.$$
It is not hard to see and in \cite{ERR} we showed it that there are a polynomial $P_l$ of degree $l+1$ and a polynomial $Q_l$ of degree $l$ such that
$$A_l(x) = \frac{P_l(x)}{(1+x^2)^{(2l+3)/2}} = \frac{(1+x^2)P_{l-1}^\prime(x) - (2l+1)x P_{l-1}(x)}{(1+x^2)^{(2l+3)/2}}$$
and
$$B_l(x) = \frac {Q_l(x)}{(1+x^2)^{(2l+3)/2}} = \frac {(1+x^2)Q_{l-1}^\prime(x) - (2l+1)x Q_{l-1}(x)}{(1+x^2)^{(2l+3)/2}}\,.$$
In \cite{ERR} we also showed that both $A_l$ and $B_l$ have only real zeros, and the zeros of $A_l$ and the zeros of $B_l$ are strictly interlacing, 
hence they are distinct. But in the present setting, we are taking linear combinations of forms that use these. So it turns out that the interlacing 
property is no longer true.

\medskip

\noindent Here is another version of \ref{Conjecture1}.

\medskip

\begin{conj} Suppose $\gamma_l$ are real numbers which are not all zero. 
The positive real zeros of the polynomials 
$$\sum_{l=0}^L{\gamma_l (1+x^2)^{L-l}x^lP_l(x)} \qquad \text{and} \qquad \sum_{l=0}^L{\gamma_l(1+x^2)^{L-l}x^lQ_l(x)}$$ 
are distinct.
\end{conj}

\begin{rem} Sometimes with forms like these one can hope for the zeros of the two expressions to interlace. Indeed, sometimes also all 
the zeros are real numbers. But in this case examples show that neither of these things are true, while it seems likely that the zeros are 
distinct (even the complex valued ones).
\end{rem}

\begin{rem} There is an additional formula that may be useful:
$$P_l(x) = lQ_{l-1}(x)(1+x^2) + xQ_l(x)\,.$$  
\end{rem}

\subsection{Conclusion: Common Zeros for $F = (X,Y)$ in a Disc}\label{CONCLUSION}

How do we now show that there can only be a finite number of zeros in a given bounded disk? In this case, it means showing there cannot be a curve of zeros anywhere. 
Indeed, if we consider $S = X^2 + Y^2$, then the zero set $S = 0$ is an analytic variety in the plane and hence, in this case, is a locally finite collection of points 
and curves. We know there are no points in this zero set outside a large disk. So all we have to do is show that there are no curves where $X$ and $Y$ are zero.

The calculation might resemble what was done in the Special Case of point charges on a line in Section 3.1, but up to now this has not worked out successfully. 
Again, this would then lead to a fairly simple complete proof of the conjecture that in the plane the electrical field from a finite number of point charges only has a 
finite number of field $F$ values that are $(0,0)$.  

On the other hand, the article by Abanov, Hayford, Khavinson, and Teodorescu~\cite{AHKT} gives a short proof in Proposition 4.1 that there are no curves in the common zero set. 
Their argument is that the following is a general fact. Take a harmonic function $\phi$ in $\mathbb R^3$, like the potential of a finite point charge electrical field.  
Assume that $\displaystyle{\frac{\partial \phi}{\partial z}}$ is zero in the $xy$-plane. Then there is no non-trivial curve in the $xy$ plane which is also simultaneously 
in both of the sets $\displaystyle{\frac{\partial \phi}{\partial x} =0}$ and $\displaystyle{\frac{\partial \phi}{\partial y} = 0}$.  
However, consider the class of functions
$$U(x,y,z) = \cosh (az) \cos (bx+cy)\,.$$

\noindent Then $U_z(x,y,0) = 0$ for all $(x,y)$ because $\sinh (0) = 0$.
Also, both $U_x$ and $U_y$ are zero for all $(x,y,z)$ such that $bx+cy = 0$. We do not even need to restrict to where $z=0$ for this to hold. 
So $U_x$ and $U_y$ are zero on the plane $\mathcal P$ where $bx+cy = 0$ and of course then the line where $bx+cy = 0$ and $z = 0$.  But importantly 
$U_z = a\sinh (az)$ on the plane $\mathcal P$ where $bx+cy = 0$.  So the gradient of $U$ does not vanish on this plane unless $a=0$.  See ~\cite{Jan}.
Moreover, $U$ is harmonic if $a^2 = b^2+c^2$.  So take examples where at least $b$ or $c$ is not zero so that $U$ is not constant. Also, an interesting point is that the 
Hessian $H(x,y,z) = U_{xx}U_{yy} - U_{xy}^2$ is identically zero everywhere in any case. So take for example $U(x,y,z) = \cosh (5z) \cos (3x - 4y)$.  
Then we have a harmonic function for which $U_z$ is zero in the $xy$-plane and $U_x = U_y = 0$ on the line where $3x - 4y = 0$ in the $xy$-plane.

Nonetheless, if we can prove a result like Proposition 4.1 in \cite{AHKT} for the potential of a finite point charge electrical field in $\mathbb R^3$, 
then we could combine this with

\medskip

\begin{conj} The zero set of a two-dimensional finite electrical field is bounded. 
\end{conj}

\medskip

\noindent The result would be 

\medskip

\begin{conj} The zero set of a two-dimensional finite electrical field consists of only finitely many points.
\end{conj}

\medskip

This is the ultimate goal in the analysis of the structure of the zero set of a finite point charge electrical field in the plane.

\subsection{Comparison with the Result Obtained by the Product Method}

An algebraic curve in the Euclidean plane is the set of the points whose coordinates are the solutions of a bivariate polynomial equation
$P(x,y) = 0$. This equation is often called the implicit equation of the curve, in contrast to the curves that are the graph of a function
defining explicitly $y$ as a function of $x$, or vice versa. With a curve given by such an implicit equation, the first problem would be to determine
the shape of the curve and to ``draw it''. These problems are not as easy to solve as in the case of the graph of a function, for which $y$ may easily
be computed for various values of $x$. The fact that the defining equation is a polynomial implies that the curve has some structural properties that
may help in solving these problems. Every algebraic curve may be uniquely decomposed into a finite number of smooth monotone arcs (called {\em branches})
sometimes connected by some points sometimes called {\em remarkable points}, and possibly a finite number of isolated points called {\em acnodes}.
A smooth monotone arc is the graph of a smooth function which is defined and monotone on an open interval of the $x$-axis or the $y$-axis. In each direction,
an arc is either unbounded (usually called an {\em infinite arc}) or has an endpoint which is either a {\em singular point} (this will be defined below) or a
point with a tangent parallel to one of the coordinate axes. The singular points of a curve of degree $d$ defined by a polynomial
$P(x,y)$ of degree $d$ are the solutions of the system of equations
$$\frac{\partial{P}}{\partial{x}}(x,y) = \frac{\partial{P}}{\partial{y}}(x,y) = P(x,y) = 0\,.$$
Harnack's Curve Theorem (see \cite{Gibson}, for instance) gives the possible numbers of connected components that an algebraic curve can have, in terms of the
degree of the curve. For any algebraic curve of degree $m$ in the real projective plane, the number of components $c$ is bounded by 
$$\frac{1-(-1)^m}{2} \leq c \leq \frac{(m-1)(m-2)}{2} + 1\,.$$

\begin{prop} Suppose  $a_1, a_2, \ldots, a_M$ are real numbers which are not all zero. Neither the zero set 
$\{X=0\}$ nor the zero set $\{Y=0\}$ can have more than $9(M-1)^24^{M-1}+1$ asymptotic directions. 
\end{prop}

These are much weaker upper bounds we have proved in the previous sections for the number of possible asymptotic directions 
for the zero sets $\{X=0\}$ and $\{Y=0\}$, but we think the ``product method" of proof below is interesting. In fact, in \cite{ERR} 
this ``product method" has been exploited in a more sophisticated way.  

\begin{proof} 
Let, as before,
$$X(x,y) = \sum_{j=1}^M{\frac{a_j(x-x_j)}{((x-x_j)^2 + y^2)^{3/2}}}\,, \quad \text{and} \quad  Y(x,y) = \sum_{j=1}^M{\frac{a_j(y-y_j)}{((x-x_j)^2 + y^2)^{3/2}}}\,.$$
Let $\Sigma$ be the collection of the $2^M$ functions $\sigma: \{1,2,\ldots,M\} \rightarrow \{-1,1\}$.
Let
$$X_{\sigma}(x,y) := \sum_{j=1}^M{\frac{\sigma(j)a_j(x-x_j)}{((x-x_j)^2 + (y-y_j)^2)^{3/2}}}\,, \quad \text{and} 
\quad Y_{\sigma}(x,y) := \sum_{j=1}^M{\frac{\sigma(j)a_j(y-y_j)}{((x-x_j)^2 + (y-y_j)^2)^{3/2}}}\,.$$
Let
$$D_j(x,y) := \prod_{k=1, k \neq j}^M{((x-x_k)^2 + (y-y_k)^2))^{3/2}}\,, \qquad j=1,2,\ldots,M\,,$$
and
$$D(x,y) := \prod_{k=1}^M{((x-x_k)^2 + (y-y_j)^2))^{3/2}}\,.$$
We have
$$X_{\sigma}(x,y) = \sum_{j=1}^M{\frac{\sigma(j)a_j(x-x_k)D_j(x,y)}{D(x,y)}}$$
and
$$Y_{\sigma}(x,y) = \sum_{j=1}^M{\frac{\sigma(j)a_j(y-y_k)D_j(x,y)}{D(x,y)}}\,.$$
Observe that the functions 
$$P(x,y) := \prod_{\sigma \in \Sigma}{\Bigg( \sum_{j=1}^M{\sigma(j)a_j(x-x_j)D_j(x,y)} \Bigg)}$$ 
and 
$$Q(x,y) := \prod_{\sigma \in \Sigma}{\Bigg( \sum_{j=1}^M{\sigma(j)a_j(y-y_j)D_j(x,y)} \Bigg)}$$
are even polynomials in each of the variables
$$\chi_j := D_j(x,y)\,, \qquad j=1,2,\ldots,M\,,$$
as they remain the same when $\chi_j$ is replaced by $-\chi_j$. Hence $X(x,y)$ and $Y(x,y)$ are polynomials in each of the variables
$$\chi_j^2 = D_j(x,y)^2, \qquad j=1,2,\ldots,M\,.$$
We conclude that
$$\prod_{\sigma \in \Sigma}{X_{\sigma}(x,y)} = \frac{P(x,y)}{D(x,y)^{2^M}} \qquad \text{and} \qquad \prod_{\sigma \in \Sigma}{Y_{\sigma}(x,y)} =  \frac{Q(x,y)}{D(x,y)^{2^M}}\,,$$
where $P$ and $Q$ are polynomials of degree at most $3(M-1)2^{M-1}$. Observe that 
$$\{X=0\} \subset \{P=0\} \qquad \text{and} \qquad \{Y=0\} \subset \{Q=0\}\,.$$  
Hence Harnack's Curve Theorem implies that neither $\{X=0\}$ nor $\{Y=0\}$ can have more than $9(M-1)^24^{M-1}+1$ asymptotic directions.
\end{proof}

\medskip

\bibliographystyle{amsplain}

\bigskip

\noindent {\bf Authors}:
\medskip

\noindent T. Erd\'elyi, Department of Mathematics, Texas A \& M University, College Station, TX,

terdelyi@math.tamu.edu
\medskip

\noindent J. Rosenblatt, Department of Mathematics, University of Illinois at Urbana-Champaign, Urbana, IL,

rosnbltt@illinois.edu
\medskip

\noindent R. Rosenblatt, Education and Human Resources at the National Science Foundation, Alexandria, VA, 

rrosenbl@nsf.gov
\bigskip

\noindent {\bf Date}: \noindent July 31, 2021

\end{document}